\documentclass[12pt]{article}

\usepackage{amsmath}
\usepackage{amsfonts}
\usepackage{amsthm}
\usepackage{hyperref}
\usepackage{amssymb}

\begin{document}

\newcommand{\A}{\mbox{${{{\cal A}}}$}}
\newcommand{\Ato}{\mbox{ ${{{\overset{\A}{\longrightarrow}}}}$ }}


\author{Attila Losonczi}
\title{Small perturbations on means and quasi-means}

\date{\today}

\newtheorem{thm}{\qquad Theorem}[section]
\newtheorem{prp}[thm]{\qquad Proposition}
\newtheorem{lem}[thm]{\qquad Lemma}
\newtheorem{cor}[thm]{\qquad Corollary}
\newtheorem{rem}[thm]{\qquad Remark}
\newtheorem{ex}[thm]{\qquad Example}
\newtheorem{df}[thm]{\qquad Definition}
\newtheorem{prb}{\qquad Problem}

\newtheorem*{thm2}{\qquad Theorem}
\newtheorem*{cor2}{\qquad Corollary}
\newtheorem*{prp2}{\qquad Proposition}

\maketitle

\begin{abstract}

\noindent

We start to investigate how small changes on the definition of ordinary means affect their properties. Especially the property of being a mean. In that direction we are looking for weakenings of the basic defining property of means. Hence we introduce weaker notions in two directions. We investigate such functions and provide many examples as well.

\noindent
\footnotetext{\noindent
AMS (2020) Subject Classifications: 26E60, 28A10, 26D99 \\

}

\end{abstract}

\section{Introduction}
In this paper we are going to extend the notion of mean in some ways.

The original motivations usually come from statistics, physics or numerical analysis. The other stream of such generalization is from when we slightly modify the defining formula of an ordinary mean. I.e. we change the function slightly and it will be no longer a mean however it is still close to being a mean in some sense.

In some of the cases if we restrict the domain of the modified function then we may end up with a new mean however on the whole domain of the original function it will not be a mean just a weaker type of mean, a quasi-mean that satisfies just one of the defining inequalities of a mean.

\smallskip

A different type of generalization arises when none of the defining inequalities hold but the function value is always close to the boundaries in some sense.

\smallskip

As a kind of application it is important to mention that there are cases when quasi-means can be used to approximate some ordinary mean.

\medskip

In the first part of the paper we define the notion of quasi-mean and provide properties that can be interesting to investigate later.

\smallskip

Then we present many examples where we point out the motivation as well. We also investigate the various properties of the presented quasi-means.

\smallskip

In the next section we investigate an operation called duality that can be interesting in its own right. However here we concentrate on how we can get quasi-means by that method.

\smallskip

Then we examine how some operations behave regarding quasi-means and how one can get a generalization of well-known theorems of means to quasi-means. For example we will examine compounding and how one can extend a $2$-variable quasi-mean to a $3$-variable quasi-mean.

\smallskip

Next we are going to deal with some generalization of monotonicity for quasi-means.

\smallskip

In the last part of the paper we provide a further possible generalization of means.

\subsection{Basic notions and notations}

An $\boldsymbol n$\textbf{-variable mean} on $H\subset\mathbb{R}$ is a function ${\cal K}:H^n\to\mathbb{R}$ such that $\forall (a_1,\dots,a_n)\in H^n \min\{a_1,\dots,a_n\}\leq {\cal K}(a_1,\dots,a_n)\leq\max\{a_1,\dots,a_n\}$ holds ($n\in\mathbb{N}-\{1\}$).
Sometimes we will use the adjective "ordinary" for such means in order to distinguish them from quasi-means. In most of the cases we may assume that the range of ${\cal K}$ is a subset of $H$ however it is not part of the definition.

\smallskip

A mean ${\cal K}$ is \textbf{reflexive} if ${\cal K}(a,\dots,a)=a$. It is \textbf{strict} if $\min\{a_1,\dots,a_n\}<{\cal K}(a_1,\dots,a_n)<\max\{a_1,\dots,a_n\}$ holds whenever $a_1,\dots,a_n$ are not all equal. It is \textbf{monotone} if $a_1\leq b_1,\dots,a_n\leq b_n$ implies that ${\cal K}(a_1,\dots,a_n)\leq{\cal K}(b_1,\dots,b_n)$. It is \textbf{continuous} if 
\[\lim\limits_{(x_1,\dots,x_n)\to(a_1,\dots,a_n)}{\cal K}(x_1,\dots,x_n)={\cal K}(a_1,\dots,a_n).\]
It is \textbf{symmetric} if $\sigma:\{1,\dots,n\}\to\{1,\dots,n\}$ being bijective implies that ${\cal K}(a_1,\dots,a_n)={\cal K}\big(a_{\sigma(1)},\dots,a_{\sigma(n)}\big)$.

An $n$-variable symmetric mean is $\boldsymbol m$\textbf{-associative} ($1\leq m\leq n$) if 
\[{\cal K}(a_1,\dots,a_n)={\cal K}\big({\cal K}(a_1,\dots,a_m),\dots,{\cal K}(a_1,\dots,a_m),a_{m+1},\dots,a_n\big)\]
where the first $m$ terms are all ${\cal K}(a_1,\dots,a_m)$.

\begin{df}
An expression ${\cal K}(a_1,\dots,a_n)$ of fixed parameters $a_1,\dots,a_n$ is \textbf{mean-like} if 
$\min\{a_1,\dots,a_n\}\leq {\cal K}(a_1,\dots,a_n)\leq\max\{a_1,\dots,a_n\}$. Here we does not consider $a_1,\dots,a_n$ as variables that can change, instead we consider them as fixed values and the expression is mean-like for those fixed values (i.e. it may not be mean-like for other values).\qed
\end{df}

A \textbf{multiset} (or mset) is a pair $\langle H,m\rangle$ where $H$ is a set and $m:H\to\mathbb{N}$. It generalizes the common notion of set in a way that it allows multiple occurrences of elements of a set (see \cite{blizard}). We say that $\langle K,l\rangle$ is a subset of $\langle H,m\rangle$ if $K\subset H$ and for each $k\in K\ \ l(k)\leq m(k)$ holds. For multisets we use the notation $\langle K,l\rangle\subset \langle H,m\rangle$.

\smallskip

If a mean ${\cal K}$ is symmetric then we can consider ${\cal K}$ as if it acted on finite multisets since the order of elements does not matter. I.e. we can write ${\cal K}(H)$ instead of ${\cal K}(a_1,\dots,a_n)$ if $H$ is the multiset of $\{a_1,\dots,a_n\}$. We will apply this convention later.

\begin{df}
Let $H\subset\mathbb{R}$ and if ${\cal K}$ is a function such that ${\cal K}$ is defined on $H^n$ for all $n\in\mathbb{N}$ (or at least for $n\geq 2$) i.e. $\text{Dom }{\cal K}=\bigcup\limits_{k=1}^{\infty}H^n$ then ${\cal K}$ is called a \textbf{multi-variable mean}. If $n\in\mathbb{N}$ then ${\cal K}|_n$ will denote ${\cal K}|_{H^n}$ i.e. it denotes its $n$-variable version.\qed
\end{df}

Throughout this paper function $\A(\dots)$ will denote the arithmetic mean of any (finite) number of variables.  

\smallskip

Let $n\in\mathbb{N}-\{1\},\ a_1,\dots,a_n\in\mathbb{R},\ f:\mathbb{R}\to\mathbb{R}$ be continuous and increasing (i.e. a homeomorphism). The quasi-arithmetic mean by $f$ is defined as 
\[{\cal K}_f(a_1,\dots,a_n)=f^{-1}\left({\frac{f(a_1)+\dots +f(a_n)}{n}}\right).\]

\smallskip

Finally we mention some usual notations: 

$\overline{\mathbb{R}}=\mathbb{R}\cup\{-\infty,\infty\},$

$\mathbb{R}^+=\{x\in\mathbb{R}:x>0\},\ \mathbb{R}^-=\{x\in\mathbb{R}:x<0\}$.

If $x\in\mathbb{R}$ then we use $|x|^+=\max\{0,x\}, |x|^-=-\min\{0,x\}$.

Set $S(x,\varepsilon)=\{y \in\mathbb{R}^n:d(x,y)<\varepsilon\}\ (x\in\mathbb{R}^n,\varepsilon>0)$.

If $x\in\mathbb{R}$ then $\lfloor x\rfloor=\max\{z\in\mathbb{Z}:z\leq x\},\ \lceil x\rceil=\min\{z\in\mathbb{Z}:x\leq z\}$.

\section{Definition and properties of quasi-means}

\begin{df}Let ${\cal K}:H^n\to\mathbb{R}\ (H\subset\mathbb{R},\ n\in\mathbb{N},\ n\geq 2)$. Then ${\cal K}$ is called a left-mean if $\min\{a_1,\dots,a_n\}\leq{\cal K}(a_1,\dots,a_n)$ holds; it is called a right-mean if ${\cal K}(a_1,\dots,a_n)\leq\max\{a_1,\dots,a_n\}$ holds ($a_1,\dots,a_n\in H$). ${\cal K}$ is called a quasi-mean if it is either a left-mean or a right-mean.
\end{df}

\begin{prp}${\cal K}$ is a mean iff it is a left-mean and a right-mean at the same time.\qed
\end{prp}

\begin{rem}In some of the cases we will assume that the range of ${\cal K}$ is a subset of $H$. Of course in those cases we will explicitly mention that.
\end{rem}

The properties reflexive, monotone, continuous, symmetric, $m$-associative are defined exactly the same way as they have been defined for ordinary means. 

\begin{df}A right-mean (left-mean) is \textbf{strict} if ${\cal K}(a_1,\dots,a_n)<\max\{a_1,\dots,a_n\}\ \ (\min\{a_1,\dots,a_n\}<{\cal K}(a_1,\dots,a_n))$ holds whenever $a_1,\dots,a_n$ are not all equal.
\end{df}

\begin{df}An $n$-variable quasi-mean ${\cal K}$ is \textbf{semi-reflexive} if $b={\cal K}(a,\dots,a)$ implies that
${\cal K}(b,\dots,b)=b$.
\end{df}

Evidently a reflexive quasi-mean is semi-reflexive.

\begin{prp}If an $n$-variable symmetric mean is $n$-associative then it is semi-reflexive.\qed
\end{prp}

\begin{df}An $n$-variable quasi-mean ${\cal K}$ is \textbf{left-continuous} if
\[\lim\limits_{\substack{x_1\to a_1, x_1\leq a_1 \\ ... \\ x_n\to a_n, x_n\leq a_n}}{\cal K}(x_1,\dots,x_n)={\cal K}(a_1,\dots,a_n).\]
Right-continuity is defined similarly with $x_i\geq a_i\ (1\leq i\leq n)$.
\end{df}

Clearly a continuous quasi-mean is both left- and right-continuous but the converse is not true.

\begin{df}A right-mean (left-mean) ${\cal K}$ is called \textbf{strong} if it satisfies that $\forall (a_1,\dots,a_n)\in\mathrm{Dom\ }{\cal K}\ \ {\cal K}(a_1,\dots,a_n)\leq\min\{a_1,\dots,a_n\}$\ \ \  ($\max\{a_1,\dots,a_n\}\leq{\cal K}(a_1,\dots,a_n)$).
\end{df}

\begin{df}A quasi-mean ${\cal K}$ is called \textbf{mean-continuous} if
\[\lim\limits_{\substack{(x_1,\dots,x_n)\to (a,\dots, a) \\ (x_1,\dots,x_n)\ne (a,\dots, a)}}{\cal K}(x_1,\dots,x_n)=a.\]
\end{df}

Obviously a mean is always mean-continuous however a quasi-mean is not. Clearly the limit is $\leq a$ for a right-mean and it is $\geq a$ for a left-mean.

\begin{prp}\label{pmeancont}A left-continuous quasi-mean is mean-continuous iff it is reflexive. (Same is true for right-continuity.)\qed
\end{prp}


\medskip

We can measure how much a quasi-mean is far from being a mean. Of course there may be many such measurements. Here we present two types.

\begin{df}Let ${\cal K}:H^n\to\mathbb{R}$ be a quasi-mean $(H\subset\mathbb{R},\ n\in\mathbb{N},\ n\geq 2)$. Set
\[\mathrm{mdist}({\cal K})=\sup\left\{\left|{\cal K}(a_1,\dots,a_n)-\max\{a_1,\dots,a_n\}\right|^+:a_1,\dots,a_n\in H\right\}+\]
\[\sup\left\{\left|\min\{a_1,\dots,a_n\}-{\cal K}(a_1,\dots,a_n)\right|^+:a_1,\dots,a_n\in H\right\},\]
\[\mathrm{mdistp}({\cal K})=\sup\left\{\frac{\left|{\cal K}(a_1,\dots,a_n)-\max\{a_1,\dots,a_n\}\right|^+}{\max\{a_1,\dots,a_n\}-\min\{a_1,\dots,a_n\}}:a_1,\dots,a_n\in H\text{ not all equal}\right\}+\]
\[\sup\left\{\frac{\left|\min\{a_1,\dots,a_n\}-{\cal K}(a_1,\dots,a_n)\right|^+}{\max\{a_1,\dots,a_n\}-\min\{a_1,\dots,a_n\}}:a_1,\dots,a_n\in H\text{ not all equal}\right\}.\]  
\end{df}

Clearly for left-means the second term is 0 while for right-means the first term is 0 for both $\mathrm{mdist}({\cal K})$ and $\mathrm{mdistp}({\cal K})$. Obviously ${\cal K}$ is a mean iff $\mathrm{mdist}({\cal K})=0$.

\begin{df}Let ${\cal K}:[a,b]^n\to\mathbb{R}$ be a quasi-mean such that the function ${\cal K}$ is Lebesgue measurable $(n\in\mathbb{N},\ n\geq 2)$. Set
\[\mathrm{mdista}({\cal K})=\frac{\lambda\big\{(a_1,\dots,a_n)\in [a,b]^n:{\cal K}(a_1,\dots,a_n)>\max\{a_1,\dots,a_n\}\big\}}{(b-a)^n}+\]
\[\frac{\lambda\big\{(a_1,\dots,a_n)\in [a,b]^n:{\cal K}(a_1,\dots,a_n)<\min\{a_1,\dots,a_n\}\big\}}{(b-a)^n}\]
where $\lambda$ is the Lebesgue measure.
\end{df}

Clearly for left-means the second term is 0 while for right-means the first term is 0 for $\mathrm{mdista}({\cal K})$. If ${\cal K}$ is a mean then $\mathrm{mdista}({\cal K})=0$. Note that the converse is not true: ${\cal K}$ can be a proper quasi-mean with $\mathrm{mdista}({\cal K})=0$. However with an obvious additional property it holds.

\begin{prp}Let ${\cal K}:[a,b]^n\to\mathbb{R}$ be a continuous quasi-mean $(n\in\mathbb{N},\ n\geq 2)$. Then if $\mathrm{mdista}({\cal K})=0$ then ${\cal K}$ is a mean.
\end{prp}
\begin{proof}Suppose indirectly that ${\cal K}$ is not a mean. Then there is $(a_1,\dots,a_n)\in [a,b]^n$ such that ${\cal K}(a_1,\dots,a_n)>\max\{a_1,\dots,a_n\}$ (the other case is similar) I.e. ${\cal K}(a_1,\dots,a_n)-\max\{a_1,\dots,a_n\}>0$. Obviously ${\cal K}-\max$ is continuous too hence there is $\delta>0$ such that $(b_1,\dots,b_n)\in S\big((a_1,\dots,a_n),\delta\big)$ implies that ${\cal K}(b_1,\dots,b_n)-\max\{b_1,\dots,b_n\}>0$.
\end{proof}

\begin{rem}If a quasi-mean is given then we can easily create an ordinary mean from it, just take 
\[\tilde{\cal K}(a_1,\dots,a_n)=\min\Big\{\max\big\{\min\{a_1,\dots,a_n\},{\cal K}(a_1,\dots,a_n)\big\},\max\{a_1,\dots,a_n\}\Big\}.\]
It truncates it on one side but keeps its behavior on the other side.\qed
\end{rem}

\section{Examples}

\subsection{Quasi-means from statistics}

Bessel's correction in statistics suggests an example for a quasi-mean. Suppose we have a statistical sample of $n$ elements, $x_1,\dots,x_n\in\mathbb{R}$. Then the sample mean is 
\[\overline{x}=\frac{\sum\limits_{i=1}^nx_i}{n}.\]
The sample standard variance is
\[\frac{\sum\limits_{i=1}^n(x_i-\overline{x})^2}{n},\]
while the unbiased sample standard variance is
\[\frac{\sum\limits_{i=1}^n(x_i-\overline{x})^2}{n-1}.\]
This suggests the following quasi-means:
\begin{ex}\label{ebessel}\[{\cal B}^+(a_1,\dots,a_n)=\frac{\sum\limits_{i=1}^na_i}{n-1}\ \ (a_i\in\mathbb{R}^+);\ {\cal B}^-(a_1,\dots,a_n)=\frac{\sum\limits_{i=1}^na_i}{n-1}\ \ (a_i\in\mathbb{R}^-).\]
Clearly $\min\{a_1,\dots,a_n\}\leq{\cal B}^+(a_1,\dots,a_n)$ and ${\cal B}^-(a_1,\dots,a_n)\leq\max\{a_1,\dots,a_n\}$ hold hence they are quasi-means while ${\cal B}^+(1,2)=3$ and ${\cal B}^-(-1,-2)=-3$ show that they are not means.
\end{ex}

\begin{prp}${\cal B}^+$ and ${\cal B}^-$ are strict, monotone, symmetric, continuous but they are not reflexive. \qed
\end{prp}

\begin{ex}There is a sequence $(a_n)$ such that 
\newline $\forall n\ {\cal B}^+(a_1,\dots,a_n)\not\leq\max\{a_1,\dots,a_n\}$. We show two such sequences.

1. Let $a_1=1,\ a_n=2\ (n>1)$. Then ${\cal B}^+(a_1,\dots,a_n)=2+\frac{1}{n-1}$.

2. Let $a_1=1,\ a_n=2-\frac{1}{2^{n-1}}\ (n>1)$. Then ${\cal B}^+(a_1,\dots,a_n)=2+\frac{1}{(n-1)2^n}$.

Similar examples can be made for ${\cal B}^-$.
\end{ex}

\begin{lem}\label{lbessel}Let $a_i\in\mathbb{R}^+,\ a_1\leq\dots\leq a_n$. Then ${\cal B}^+(a_1,\dots,a_n)$ is mean-like iff ${\cal B}^+(a_1,\dots,a_{n-1})\leq a_n$.
\end{lem}
\begin{proof}${\cal B}^+(a_1,\dots,a_n)\leq a_n$ iff $\sum\limits_{i=1}^na_i\leq(n-1)a_n$ iff $\sum\limits_{i=1}^{n-1}a_i\leq(n-2)a_n$ iff ${\cal B}^+(a_1,\dots,a_{n-1})\leq a_n$.
\end{proof}

\begin{prp}Let $a_i\in\mathbb{R},\ a_1\leq\dots\leq a_n$. If $a_1\leq 0,a_n\geq 0$ then 
\[a_1\leq\frac{\sum\limits_{i=1}^na_i}{n-1}\leq a_n.\]
\end{prp}
\begin{proof}\[a_1\leq\frac{\sum\limits_{i=1}^{n-1}a_i}{n-1}\leq\frac{\sum\limits_{i=1}^na_i}{n-1}\leq\frac{\sum\limits_{i=2}^na_i}{n-1}\leq a_n.\qedhere\]
\end{proof}

\begin{thm}\label{tbpati}If $a_n\in\mathbb{R}^+,\ a_n\to+\infty$ then $\exists n\in\mathbb{N}$ such that ${\cal B}^+(a_1,\dots,a_n)$ is mean-like.
\end{thm}
\begin{proof}Suppose the contrary: $\forall n\in\mathbb{N}\ \max\{a_1,\dots,a_n\}<\frac{\sum\limits_{i=1}^{n}a_i}{n-1}$. Take such $n$. Obviously
\[a_n<\frac{\sum\limits_{i=1}^{n}a_i}{n-1}\]
is equivalent to
\[a_n<\frac{\sum\limits_{i=1}^{n-1}a_i}{n-2}.\]
Now we get that
\[(n-2)a_n<\sum\limits_{i=1}^{n-1}a_i\]
\[(n-2)\sum\limits_{i=1}^{n}a_i=(n-2)\sum\limits_{i=1}^{n-1}a_i+(n-2)a_n<(n-1)\sum\limits_{i=1}^{n-1}a_i\]
\[\frac{\sum\limits_{i=1}^{n}a_i}{n-1}<\frac{\sum\limits_{i=1}^{n-1}a_i}{n-2}\]
but by the indirect assumption  
\[a_{n+1}<\frac{\sum\limits_{i=1}^{n-1}a_i}{n-2}\]
hence $\forall m\geq n$
\[a_{m}<\frac{\sum\limits_{i=1}^{n-1}a_i}{n-2}\]
which is a contradiction.
\end{proof}

\begin{thm}\label{tbptail}If ${\cal B}^+(a_1,\dots,a_n)\leq\max\{a_1,\dots,a_n\}$ then ${\cal B}^+(a_1,\dots,a_n,a_{n+1})\leq\max\{a_1,\dots,a_n,a_{n+1}\}\ \ (a_i\in\mathbb{R}^+)$.
\end{thm}
\begin{proof}Let $a_1\leq\dots\leq a_n$.
Then $\frac{\sum\limits_{i=1}^{n}a_i}{n-1}\leq a_n$ holds. There are two cases:
\begin{enumerate}\setlength\itemsep{0em}
\item If $a_n\leq a_{n+1}$ then $\frac{\sum\limits_{i=1}^{n}a_i}{n-1}\leq a_n\leq a_{n+1}$ gives that $\sum\limits_{i=1}^{n}a_i\leq (n-1)a_{n+1}$ and $\sum\limits_{i=1}^{n+1}a_i\leq na_{n+1}$ which yields that $\frac{\sum\limits_{i=1}^{n+1}a_i}{n}\leq a_{n+1}$. 
\item If $a_{n+1} < a_n$ then $\max\{a_1,\dots,a_n,a_{n+1}\}=a_n$. Now summing $\sum\limits_{i=1}^{n}a_i\leq (n-1)a_{n}$ and $a_{n+1} < a_n$ gives that $\sum\limits_{i=1}^{n+1}a_i\leq na_{n}$ hence $\frac{\sum\limits_{i=1}^{n+1}a_i}{n}\leq a_n$ holds. \qedhere
\end{enumerate}
\end{proof}

\begin{cor}If $a_n\in\mathbb{R}^+,\ a_n\to+\infty$ then $\exists N\in\mathbb{N}$ such that $n<N$ implies that ${\cal B}^+(a_1,\dots,a_n)$ is not mean-like while $n\geq N$ implies that ${\cal B}^+(a_1,\dots,a_n)$ is mean-like.
\end{cor}
\begin{proof}\ref{tbpati} and \ref{tbptail}.
\end{proof}

When we consider the unbiased sample standard deviation
\[\sqrt{\frac{\sum\limits_{i=1}^n(x_i-\overline{x})^2}{n-1}}\]
then we can also produce the following quasi-mean:
\begin{ex}\[{\cal K}(a_1,\dots,a_n)=\sqrt{\frac{\sum\limits_{i=1}^na_i^2}{n-1}}\ \ (a_i\in\mathbb{R}^+).\]
\end{ex}

We will investigate this kind of quasi-means in section \ref{sqaqmeans} in a more general aspect.

\medskip

In statistics sometimes we get rid of the smallest and/or the largest elements of a sample in order to get a more usable sample. This motivates the definition of the next quasi-means. We restrict ourselves for the positive cases only.

\begin{ex}Let $a_1,\dots,a_n\in\mathbb{R}^+,\ a_1\leq\dots\leq a_n,\ n\geq 3$. Then set
\[{\cal K}_1(a_1,\dots,a_n)=\frac{\sum\limits_{i=2}^na_i}{n},\]
\[{\cal K}_2(a_1,\dots,a_n)=\frac{\sum\limits_{i=1}^{n-1}a_i}{n},\]
\[{\cal K}_3(a_1,\dots,a_n)=\frac{\sum\limits_{i=2}^{n-1}a_i}{n}.\]
Obviously ${\cal K}_3(a_1,\dots,a_n)\leq{\cal K}_2(a_1,\dots,a_n)\leq{\cal K}_1(a_1,\dots,a_n)\leq\max\{a_1,\dots,a_n\}$ hence all of them are quasi-means. The example ${\cal K}_1(10,11,12)<10$ demonstrates that they are not means.
\end{ex}

\subsection{Quasi-means from numerical analysis}

The following quasi-mean inspired by approximation of ordinary means i.e. the motivation comes from numerical analysis. 

Suppose we have $n$ real numbers, $a_1,\dots,a_n\in\mathbb{R}$ and we want to calculate their arithmetic mean e.g. But among the numbers there are some irrational numbers too therefore the calculation is possible only approximately. In order to do so we take approximations for all numbers, e.g. we take the first three digits of each number and use them in the calculation. The more digits we take, the better approximation we get.

\begin{ex}\label{eqmna}Let $a_i\in\mathbb{R}\ (1\leq i\leq n)$ and $a_i=x_i.x_{i,1}x_{i,2}x_{i,3}\dots$ where $x_i\in\mathbb{Z}, x_{i,k}\in\{0,\dots,9\}\ (k\in\mathbb{N})$ i.e. we take the decimal representation of each number where we also assume that $\forall i\ \forall N\ \exists k>N$ such that $x_{i,k}\ne 9$ i.e. the representation is without trailing 9's. 
Then set
\[{\cal A}_m(a_1,\dots,a_n)=\A(x_1.x_{1,1}\dots x_{1,m};\dots;x_n.x_{n,1}\dots x_{n,m})\ \ (m\in\mathbb{N}\cup\{0\})\]
in other words we take the first $m$ digits of each number and calculate the arithmetic mean of those numbers. In that way we get an approximation of $\A(a_1,\dots,a_n)$. 

Let us emphasize ${\cal A}_0$ that is ${\cal A}_0(a_1,\dots,a_n)=\A(\lfloor a_1\rfloor ,\dots, \lfloor a_n\rfloor )$ i.e. we take the arithmetic mean of the floors of the parameter values.
\end{ex}

Here we provide a more precise equivalent definition that allows some generalization as well.

\begin{df}\label{dqmna}Let $m\in\mathbb{Z}$. Let 
\[\text{Dom }\A_m=\{(a_1,\dots,a_n)\in\mathbb{R}^n\ (n\in\mathbb{N}): \exists i\in\{1,\dots,n\}\text{ such that }\lfloor10^ma_i\rfloor\ne 0\}.\]
If $(a_1,\dots,a_n)\in\text{Dom }\A_m$ then set
\[\A_m(a_1,\dots,a_n)=\A\left(\frac{\lfloor 10^ma_1\rfloor}{10^m},\dots,\frac{\lfloor 10^ma_n\rfloor}{10^m}\right).\]
We will also use the notation ${\cal A}^-_m$ for ${\cal A}_m$.\qed
\end{df}

\begin{rem}1) Notice that we excluded those tuples $(a_1,\dots,a_n)$ where $\forall i\ \lfloor10^ma_i\rfloor=0$ since for such tuple $\A_m(a_1,\dots,a_n)$ would be 0 that would cause difficulties and not too nice behavior.

2) Especially $\A_m$ is not defined for $(0,\dots,0)$.

3) In the definition of $\A_m$ inside \A\ we do not ignore zeros i.e. we do not drop zeros from the tuple inside \A.

4) $\A_m$ may be defined for a single number however its result can be different, e.g. $\A_0(1.1)=1$.
\end{rem}

\begin{rem}Now we have $\A_m$ defined for negative $m$ as well that still can have practical value. Because when we have to calculate the arithmetic mean of large numbers then sometimes we just care about the first 3-4 digits and not the remaining ones. Moreover if there are small values among the numbers as well then sometimes we may ignore those too, better to say we replace them with 0. These two operations altogether act like $\A_m$ for some negative $m$. 
\end{rem}

\begin{prp}For $\text{Dom }\A_m$ the following statements hold:
\begin{enumerate}\setlength\itemsep{0em}
\item $\text{Dom }\A_m=\bigcup\limits_{n=1}^{\infty}\mathbb{R}^n-[0,10^{-m})^n$
\item If $k<m$ then $\text{Dom }\A_k\subsetneq\text{Dom }\A_m$.
\item $\bigcup\limits_{m=-\infty}^{+\infty}\text{Dom }\A_m=\bigcup\limits_{n=1}^{\infty}\mathbb{R}^n-\{0\}^n$.
\item If $n\in\mathbb{N}$ then $\bigcup\limits_{m=1}^{-\infty}\mathbb{R}^n-\text{Dom }\A_m=\left(\mathbb{R}^+\cup\{0\}\right)^n$.
\end{enumerate}
\end{prp}
\begin{proof}
\begin{enumerate}\setlength\itemsep{0em}
\item $(a_1,\dots,a_n)\not\in\text{Dom }\A_m\iff\forall i\in\{1,\dots,n\}\ a_i\in[0,10^{-m})$.
\item By (1) $\text{Dom }\A_k\cap\mathbb{R}^n\subsetneq\text{Dom }\A_m\cap\mathbb{R}^n$ holds for all $n\in\mathbb{R}$.
\item Consequence of (1).
\item Consequence of (1).\qedhere
\end{enumerate}
\end{proof}

\begin{prp}If $m\in\mathbb{Z}, n\in\mathbb{N}, (a_1,\dots,a_n)\in\text{Dom }\A_0$ then 
\[\A_m(a_1,\dots,a_n)=\frac{1}{10^m}\A_0(10^ma_1,\dots,10^ma_n).\tag*{\qed}\]
\end{prp}

\begin{prp}\label{pqmna1}Let ${\cal A}_m$ be the function defined in \ref{dqmna} $(m\in\mathbb{Z})$. Then the following statements hold:
\begin{enumerate}\setlength\itemsep{0em}
\item ${\cal A}_m(a_1,\dots,a_n)\leq\max\{a_1,\dots,a_n\}$ i.e. ${\cal A}_m$ is a right-mean ($(a_1,\dots,a_n)\in\text{Dom }\A_m$)

\item ${\cal A}_m(19\cdot 10^{-m-1},21\cdot 10^{-m-1})=15\cdot 10^{-m-1}$ demonstrates that ${\cal A}_m$ is not a mean

\item ${\cal A}_m\leq{\cal A}_{m+1}$

\item ${\cal A}_m\leq{\cal A}$

\item ${\cal A}-{\cal A}_m\leq\frac{1}{10^m}$

\item $\lim\limits_{m\to\infty}{\cal A}_m(a_1,\dots,a_n)={\cal A}(a_1,\dots,a_n)$ moreover the convergence is uniform ($(a_1,\dots,a_n)\in\text{Dom }\A_0$)

\item a) $\mathrm{mdist}({\cal A}_m)=\frac{1}{10^m}$

b) $\mathrm{mdistp}({\cal A}_m)=\infty$

c) for ${\cal A}_m|_2$\ \ $\mathrm{mdista}({\cal A}_m|_{[1,2]^2})=\frac{4}{10^m}\left(1-\frac{1}{10^m}\right)$.
\end{enumerate}
\end{prp}
\begin{proof}1,3 and 4 comes from the monotonicity of \A. 

5 follows from the shift invariance of \A. 

6 is a consequence of 5.

\smallskip

To see 7.a observe that 
\[\min\{a_1,\dots,a_n\}-{\cal A}_m(a_1,\dots,a_n)\leq {\cal A}(a_1,\dots,a_n)-{\cal A}_m(a_1,\dots,a_n)\leq\frac{1}{10^m}\]
and for $a=\frac{1}{10^m}-\frac{1}{10^k}\ (k>m)$ we get that 
\[\min\{a,\dots,a\}-{\cal A}_m(a,\dots,a)=\frac{1}{10^m}-\frac{1}{10^k}\to\frac{1}{10^m}\text{ if }k\to\infty.\]

\smallskip

7.b: Let $a_1=\dots=a_{n-1}=1-\frac{1}{10^{m+1}}$ and let $a_{n_k}\to 1-\frac{1}{10^{m+1}},\ a_{n_k}>1-\frac{1}{10^{m+1}}$. Then 
\[\lim\limits_{k\to\infty}\left|\min\{a_1,\dots,a_{n-1},a_{n_k}\}-{\cal A}_m(a_1,\dots,a_{n-1},a_{n_k})\right|^+=\frac{1}{10^{m}}-\frac{1}{10^{m+1}}\]
while $\lim\limits_{k\to\infty}\max\{a_1,\dots,a_{n-1},a_{n_k}\}-\min\{a_1,\dots,a_{n-1},a_{n_k}\}=0$
which gives that $\mathrm{mdistp}({\cal A}_m)=\infty$.

\smallskip

7.c: Clearly ${\cal A}_m(x,y)\geq\min\{x,y\}$ holds iff $|x-y|\geq\frac{2}{10^{m}}$. Then 
\[\lambda\big\{(x,y)\in [1,2]^2:{\cal A}_m(x,y)\geq\min\{x,y\}\big\}=\left(1-\frac{2}{10^{m}}\right)^2\]
which gives that 
\[\lambda\big\{(x,y)\in [1,2]^2:{\cal A}_m(x,y)<\min\{x,y\}\big\}=\frac{4}{10^m}\left(1-\frac{1}{10^m}\right).\qedhere\]
\end{proof}

\begin{prp}\label{pamlim}Let $a_n\to\alpha\in\overline{\mathbb{R}}$. Then 
\begin{enumerate}\setlength\itemsep{0em}
\item If $m\in\mathbb{Z},\ \alpha\in\mathbb{R},\ 10^m\alpha\not\in\mathbb{Z}$ then $\lim\limits_{n\to\infty}\A_m(a_1,\dots,a_n)=\A_m(\alpha)$.
\item If $m\in\mathbb{Z}\ \alpha\in\mathbb{R},\ \alpha\ne 10^{-m},\ \forall n\ \alpha\leq a_n$ then 
\newline $\lim\limits_{n\to\infty}\A_m(a_1,\dots,a_n)=\A_m(\alpha)$.
\item If $\alpha=\pm\infty$ then $\lim\limits_{n\to\infty}\A_m(a_1,\dots,a_n)=\alpha\ \ (m\in\mathbb{Z})$.
\item $\lim\limits_{\substack{m\to\infty \\ n\to\infty}}\A_m(a_1,\dots,a_n)=\alpha$.
\end{enumerate}
\end{prp}
\begin{proof} We will refer to the well-known result that $a_n\to\alpha\in\overline{\mathbb{R}}$ implies that $\lim\limits_{m\to\infty}\A(a_1,\dots,a_m)=\alpha$.
\begin{enumerate}\setlength\itemsep{0em}
\item By definition
$\A_m(a_1,\dots,a_n)=\A\left(\frac{\lfloor 10^ma_1\rfloor}{10^m},\dots,\frac{\lfloor 10^ma_n\rfloor}{10^m}\right)$. Clearly $10^ma_n\to 10^m\alpha$ and then $\lfloor 10^ma_n\rfloor\to\lfloor 10^m\alpha\rfloor$ which yields that $\frac{\lfloor 10^ma_n\rfloor}{10^m}\to\frac{\lfloor 10^m\alpha\rfloor}{10^m}$ which gives the claim.
\item $\frac{\lfloor 10^ma_n\rfloor}{10^m}\to\frac{\lfloor 10^m\alpha\rfloor}{10^m}$ holds in this case too.
\item $\frac{\lfloor 10^ma_n\rfloor}{10^m}\to\pm\infty$.
\item By \ref{pqmna1} (4) and (5) $\A-\frac{1}{10^m}\leq\A_m\leq\A$.\qedhere
\end{enumerate}
\end{proof}

\begin{ex}If $10^m\alpha\in\mathbb{Z}$ and we omit the condition that $\alpha\leq a_n$ then \ref{pamlim} (2) does not hold. 

\smallskip

E.g. let $a_n=2-\frac{1}{n}$. Then $\A_0(a_1,\dots,a_n)=1\to 1\ne\A_0(2)=2$.
\end{ex}

\begin{prp}\label{ppropofAm}${\cal A}_m$ is strict, monotone, symmetric, right-continuous and semi-reflexive. It is not reflexive, not mean-continuous and not left-continuous hence not continuous either. The continuity points equals to the left-continuity points which are $\{(a_1,\dots,a_n): \forall i\ 10^{-m}a_i\not\in\mathbb{Z}\}$.
\end{prp}
\begin{proof}${\cal A}_m$ is strict because ${\cal A}_m\leq{\cal A}<\max$.

\smallskip

${\cal A}_m$ is monotone because $f(x)=\lfloor x\rfloor$ is a monotone function and \A\ is a monotone mean.

\smallskip

${\cal A}_m$ is symmetric because so is \A.

\smallskip

${\cal A}_m$ is right-continuous because $f(x)=\lfloor x\rfloor$ is a right-continuous function and \A\ is a continuous mean.

\smallskip

${\cal A}_m$ is semi-reflexive because if $b={\cal A}_m(a,\dots,a)=\frac{\lfloor10^m a\rfloor}{10^m}$ then ${\cal A}_m(b,\dots,b)=\frac{1}{10^m}\A(\lfloor10^m a\rfloor,\dots,\lfloor10^m a\rfloor)=b$.

\smallskip

${\cal A}_m$ is not reflexive as e.g. ${\cal A}_m(2.1\cdot 10^{-m},\dots,2.1\cdot 10^{-m})=2\cdot 10^{-m}$ shows.

\smallskip

${\cal A}_m$ is not mean-continuous by \ref{pmeancont} and failing to be reflexive.

\smallskip

${\cal A}_m$ is not left-continuous as take e.g. the sequence $a_k=2\cdot 10^{-m}-\frac{1}{k}$ and $\lim\limits_{k\to\infty}{\cal A}_m|_n(a_k,\dots,a_k)=10^{-m}$ while ${\cal A}_m|_n(2\cdot 10^{-m},\dots,2\cdot 10^{-m})=2\cdot 10^{-m}$.

\smallskip

Regarding (left-)continuity points if for any $i\ 10^{-m}a_i\in\mathbb{Z}$ then ${\cal A}_m|_n$ is not left-continuous at $(a_1,\dots,a_n)$ since $\lim\limits_{k\to\infty}{\cal A}_m|_n(a_1,\dots,a_{i-1},a_i-\frac{1}{k},a_{i+1},\dots,a_n)\ne{\cal A}_m|_n(a_1,\dots,a_n)$. 

If $\forall i\ 10^{-m}a_i\not\in\mathbb{Z}$ and $\forall i\ \lim\limits_{k\to\infty}a_k^{(i)}=a_i$ then 
\[\lim\limits_{k\to\infty}\frac{\lfloor 10^m a_k^{(i)}\rfloor}{10^m}=\frac{\lfloor10^m a_i\rfloor}{10^m}\]
which altogether with the fact that \A\ is continuous gives that ${\cal A}_m|_n$ is continuous at such points.
\end{proof}

We examine $k$-associativity.

\begin{prp}Let $1\leq k\leq n$. Then
\[{\cal A}_m|_n\big({\cal A}_m|_k(a_1,\dots,a_k),\dots,{\cal A}_m|_k(a_1,\dots,a_k),a_{k+1},\dots,a_n\big)\leq{\cal A}_m|_n(a_1,\dots,a_n).\]
\end{prp}
\begin{proof}Clearly
\[{\cal A}_m\big({\cal A}_m(a_1,\dots,a_k),\dots,{\cal A}_m(a_1,\dots,a_k),a_{k+1},\dots,a_n\big)=\]
\[{\cal A}\left(\frac{\left\lfloor\frac{\lfloor 10^ma_1\rfloor,\dots,\lfloor 10^ma_k\rfloor}{k}\right\rfloor,\dots,\left\lfloor\frac{\lfloor 10^ma_1\rfloor,\dots,\lfloor 10^ma_k\rfloor}{k}\right\rfloor,\lfloor 10^ma_{k+1}\rfloor,\dots,\lfloor 10^ma_n\rfloor}{10^m n}\right)\leq\]
\[{\cal A}\left(\frac{\frac{\lfloor 10^ma_1\rfloor,\dots,\lfloor 10^ma_k\rfloor}{k},\dots,\frac{\lfloor 10^ma_1\rfloor,\dots,\lfloor 10^ma_k\rfloor}{k},\lfloor 10^ma_{k+1}\rfloor,\dots,\lfloor 10^ma_n\rfloor}{10^m n}\right)=\]
\[{\cal A}\left(\frac{\lfloor 10^ma_1\rfloor,\dots,\lfloor 10^ma_k\rfloor,\lfloor 10^ma_{k+1}\rfloor,\dots,\lfloor 10^ma_n\rfloor}{10^m n}\right)={\cal A}_m(a_1,\dots,a_n).\qedhere\]

\end{proof}

\begin{ex}${\cal A}_m|_n$ is not $k$-associative in general ($1\leq k\leq n$). Set $v=10^{-m}$.

Let \[a_i=
\begin{cases}
2v&\text{if }i=k\\
nv&\text{if }i=n\\
v&\text{otherwise.}\\
\end{cases}
\]
Now ${\cal A}_m|_n(a_1,\dots,a_n)={\cal A}_m(v,\dots,v,2v,v,\dots,v,nv)=\frac{n-2+2+n}{n}v=2v$. While ${\cal A}_m|_k(a_1,\dots,a_k)={\cal A}_m|_k(v,\dots,v,2v)<2v$ and $\lfloor{\cal A}_m|_k(v,\dots,v,2v)\rfloor=v$ therefore
\[{\cal A}_m|_n({\cal A}_m|_k(v,\dots,v,2v),\dots,{\cal A}_m|_k(v,\dots,v,2v),v,\dots,v,nv)=\]
\[\A(v,\dots,v,nv)<2v.\tag*{\qed}\]
\end{ex}

\begin{prp}$\mathrm{Ran\ }{\cal A}_m|_n$ is a closed and discrete set moreover $\mathrm{Ran\ }{\cal A}_m|_n=\left\{\frac{1}{10^m}\frac{k}{n}:k\in\mathbb{Z}\right\}$.
\end{prp}
\begin{proof}First note that $\mathrm{Ran\ }{\cal A}_m|_n=\frac{1}{10^m}\mathrm{Ran\ }{\cal A}_0|_n$ because
\[\A_m(a_1,\dots,a_n)=\A\left(\frac{\lfloor 10^ma_1\rfloor}{10^m},\dots,\frac{\lfloor 10^ma_n\rfloor}{10^m}\right)=\]
\[\frac{1}{10^m}\A\left(\lfloor 10^ma_1\rfloor,\dots,\lfloor 10^ma_n\rfloor\right)=\frac{1}{10^m}\A_0\left(10^ma_1,\dots,10^ma_n\right).\]
Evidently $\mathrm{Ran\ }{\cal A}_0|_n\subset\left\{\frac{k}{n}:k\in\mathbb{Z}\right\}$ because $\A_0(a_1,\dots,a_n)=\A(\lfloor a_1\rfloor,\dots,\lfloor a_n\rfloor)$. If $k\in\mathbb{Z}$ then $\A_0|_n(1,\dots,1,k-n+1)=\frac{k}{n}$ which gives the other inclusion.
\end{proof}

Now we provide a sufficient bound for ${\cal A}_m(a_1,\dots,a_n)$ to be mean-like.

\begin{prp}Let $a_1,\dots,a_n\in\mathbb{R}$ be fixed, $a_1\leq\dots\leq a_n,\ a_1<a_n$. Then $m>\lg\left(\frac{n}{a_n-a_1}\right)+1$ implies that ${\cal A}_m(a_1,\dots,a_n)$ is mean-like.
\end{prp}
\begin{proof}First note that 
\[\frac{a_n-a_1}{n}=\frac{(n-1)a_1+a_n}{n}-a_1=\A(a_1,\dots,a_1,a_n)-a_1\]
where there are $n-1$ pieces of $a_1$ inside \A. 

Let $m>\lg(\frac{a_n-a_1}{n})+1$. Then it is clear that
\[\A(a_1,\dots,a_{n-1},a_n)-\frac{1}{10^m}\leq\A_m(a_1,\dots,a_{n-1},a_n)\leq\A(a_1,\dots,a_{n-1},a_n)\]
but
\[\A(a_1,\dots,a_1,a_n)-\frac{1}{10^m}\leq\A(a_1,\dots,a_{n-1},a_n)-\frac{1}{10^m}.\]
Now
\[\frac{1}{10^m}<\A(a_1,\dots,a_1,a_n)-a_1\]
equivalent with $a_1<\A(a_1,\dots,a_1,a_n)-\frac{1}{10^m}$ which completes the proof.
\end{proof}

\begin{cor}Let $n\in\mathbb{N}-\{1\}$ be fixed, $h>0$. Let 
\[L=\{(a_1,\dots,a_n)\in \mathbb{R}^n: \max\{a_1,\dots,a_n\}-\min\{a_1,\dots,a_n\}\geq h\}.\]
Set $\hat{\cal A}_m={\cal A}_m|_L\ (m\in\mathbb{N}\cup\{0\})$. Then there is $M\in\mathbb{N}$ such that $m>M$ implies that $\hat{\cal A}_m$ is a mean. \qed
\end{cor}

\begin{ex}In the definition \ref{dqmna} we could replace "arithmetic mean" with any continuous (monotone) ordinary mean.

E.g. for the geometric mean ${\cal G}$ we can have the following.
Let $m\in\mathbb{Z}$. Let 
\[\text{Dom }{\cal G}_m=\{(a_1,\dots,a_n)\in(\mathbb{R}^+)^n\ (n\in\mathbb{N}): \forall i\in\{1,\dots,n\}\text{ such that }\lfloor10^ma_i\rfloor\ne 0\}.\]
If $(a_1,\dots,a_n)\in\text{Dom }{\cal G}_m$ then set
${\cal G}_m(a_1,\dots,a_n)={\cal G}\left(\frac{\lfloor 10^ma_1\rfloor}{10^m},\dots,\frac{\lfloor 10^ma_n\rfloor}{10^m}\right).$

Notice the difference in the domain: $\forall i$ instead of $\exists i$. The following example shows why that is needed. 

Let us observe $(100,1)$: ${\cal G}(100,1)=10$ and if we allowed $(100,1)$ for ${\cal G}_{-1}$ then we would get that ${\cal G}_{-1}(100,1)=0$ that would bring strange behavior that we want to avoid.
\end{ex}

Now we define some counterparts of ${\cal A}_m$.

\begin{df}\label{dqmnaplus}Let $m\in\mathbb{Z}$. Let 
\[\text{Dom }\A^+_m=\{(a_1,\dots,a_n)\in\mathbb{R}^n\ (n\in\mathbb{N}): \exists i\in\{1,\dots,n\}\text{ such that }\lceil10^ma_i\rceil\ne 0\}.\]
If $(a_1,\dots,a_n)\in\text{Dom }\A^+_m$ then set
\[\A^+_m(a_1,\dots,a_n)=\A\left(\frac{\lceil 10^ma_1\rceil}{10^m},\dots,\frac{\lceil 10^ma_n\rceil}{10^m}\right).\]
\end{df}

One can prove similar statements for $\A^+_m$ as we get in \ref{pqmna1}. We omit the proof as it is similar.

\begin{prp}\label{pqmna12}For $\A^+_m$ the following statements hold:
\begin{enumerate}\setlength\itemsep{0em}
\item $\min\leq{\cal A}^+_m$ i.e. ${\cal A}^+_m$ is a left-mean

\item ${\cal A}^+_m$ is not a mean

\item ${\cal A}^+_{m+1}\leq{\cal A}^+_m$

\item ${\cal A}\leq{\cal A}^+_m$

\item ${\cal A}^+_m-{\cal A}\leq\frac{1}{10^m}$

\item $\lim\limits_{m\to\infty}{\cal A}^+_m(a_1,\dots,a_n)={\cal A}(a_1,\dots,a_n)$ moreover the convergence is uniform ($(a_1,\dots,a_n)\in\text{Dom }\A_0^+$)

\item a) $\mathrm{mdist}({\cal A}^+_m)=\frac{1}{10^m}$

b) $\mathrm{mdistp}({\cal A}^+_m)=\infty$

c) for 2-variable ${\cal A}^+_m$ $\mathrm{mdista}({\cal A}^+_m|_{[1,2]^2})=\frac{4}{10^m}\left(1-\frac{1}{10^m}\right)$.
\end{enumerate}
\end{prp}

\begin{df}Set $\A^{-+}_m=\A^-_m+\frac{1}{10^m},\ \A^{+-}_m=\A^+_m-\frac{1}{10^m}$.
\end{df}

\begin{rem}Actually $\A^{-+}_m$ corresponds to the algorithm that we truncate the numbers at the m$^{th}$ digit and increase the m$^{th}$ digit by 1 and calculate the arithmetic mean of those numbers. 

For $\A^{+-}_m$ similarly with decreasing the m$^{th}$ digit by 1.
\end{rem}

\begin{prp}$\A^{-}_m$ and $\A^{+-}_m$ are right-means. 
\newline $\A^{+}_m$ and $\A^{-+}_m$ are left-means.
\end{prp}
\begin{proof}By \ref{pqmna1} (5) $\min\leq\A\leq\A^{-+}_m$. By \ref{pqmna12} (5) $\A^{+-}_m\leq\A\leq\max$.
\end{proof}

\begin{prp}$\A^{+-}_m=10^{-m}\A^{+-}_0,\ \A^{-+}_m=10^{-m}\A^{-+}_0$.\qed
\end{prp}

\begin{prp}\label{pineqA}$\A-\frac{1}{10^m}\leq\A^{+-}_m\leq\A^{-}_m\leq \A\leq \A^{+}_m\leq\A^{-+}_m\leq\A+\frac{1}{10^m}$.
\end{prp}
\begin{proof}Clearly $\lceil x\rceil-1\leq\lfloor x\rfloor\leq\lceil x\rceil\leq\lfloor x\rfloor+1$ hence
\[\A\left(\frac{ 10^ma_1}{10^m},\dots,\frac{ 10^ma_n}{10^m}\right)-\frac{1}{10^m}\leq
\A\left(\frac{\lceil 10^ma_1\rceil-1}{10^m},\dots,\frac{\lceil 10^ma_n\rceil-1}{10^m}\right)\leq\]
\[\A\left(\frac{\lfloor 10^ma_1\rfloor}{10^m},\dots,\frac{\lfloor 10^ma_n\rfloor}{10^m}\right)\leq
\A\left(\frac{ 10^ma_1}{10^m},\dots,\frac{ 10^ma_n}{10^m}\right)\leq\]
\[\A\left(\frac{\lceil 10^ma_1\rceil}{10^m},\dots,\frac{\lceil 10^ma_n\rceil}{10^m}\right)\leq
\A\left(\frac{\lfloor 10^ma_1\rfloor+1}{10^m},\dots,\frac{\lfloor 10^ma_n\rfloor+1}{10^m}\right).\qedhere\]
\end{proof}

\begin{ex}We cannot expect equalities in general in \ref{pineqA}.

E.g. $\A^{+-}_0(2.1,3)<\A^{-}_0(2.1,3)<\A(2.1,3)< \A^{+}_0(2.1,3)<\A^{-+}_0(2.1,3)$ since 
$\A^{+-}_0(2.1,3)=2,\ \A^{-}_0(2.1,3)=2.5,\ \A(2.1,3)=2.55,\ \A^{+}_0(2.1,3)=3,\ \A^{-+}_0(2.1,3)=3.5$.

Multiplying the parameters by $10^{-m}$ provides example for $m$ instead of 0.
\end{ex}

\begin{prp}${\cal A}^{-+}_m,{\cal A}^{+-}_m$ and ${\cal A}^{+}_m$ are strict, monotone, symmetric. They are not reflexive, not mean-continuous and not continuous either. 

${\cal A}^{-+}_m$ is right-continuous, ${\cal A}^{+-}_m$ and ${\cal A}^{+}_m$ are left-continuous. ${\cal A}^{-+}_m$ is not left-continuous, ${\cal A}^{+-}_m$ and ${\cal A}^{+}_m$ are not right-continuous. ${\cal A}^{+}_m$ is semi-reflexive while ${\cal A}^{-+}_m$ and ${\cal A}^{+-}_m$ are not.
\end{prp}
\begin{proof}Everything can be proved similarly as in \ref{ppropofAm}. Only non-semi-reflexivity is new: ${\cal A}^{-+}_m(2.1\cdot 10^{-m},\dots,2.1\cdot 10^{-m})=3\cdot 10^{-m}$ and ${\cal A}^{-+}_m(3\cdot 10^{-m},\dots,3\cdot 10^{-m})=4\cdot 10^{-m}$. Similar example can be found for ${\cal A}^{+-}_m$.
\end{proof}

Let us generalize the method (def \ref{dqmna}) even further.

\begin{prp}\label{pqmna2}Let ${\cal K}$ be a continuous mean of $n$ variables. For each $x\in\mathbb{R}$ we fix a sequence $(x_i)$ such that $x_i\leq x$ and $\lim\limits_{i\to\infty}x_i=x$.
Let ${\cal K}_m(a_1,\dots,a_n)={\cal K}(a_{1,m},\dots,a_{n,m})\ (m\in\mathbb{N},\ a_1,\dots,a_n\in\mathbb{R})$ where the sequence $(a_{i,m})$ belongs to $a_i\ (1\leq i\leq n)$. 
Then ${\cal K}_m$ is a quasi-mean.
\end{prp}
\begin{proof}Let $a_{j,m}=\max\{a_{1,m},\dots,a_{n,m}\}$. Then $a_{j,m}\leq a_j\leq\max\{a_{1},\dots,a_{n}\}$ which gives that ${\cal K}_m(a_1,\dots,a_n)\leq a_{j,m}\leq \max\{a_{1},\dots,a_{n}\}$.
\end{proof}

\begin{prp}\label{pqmna3}Let $a_1,\dots,a_n\in\mathbb{R}$ be fixed and ${\cal K}$ be a strict, continuous mean. For each $i\in\{1,\dots,n\}$ let $(a_{i,m})$ be a sequence such that $\lim\limits_{m\to\infty}a_{i,m}=a_i\ (i\in\{1,\dots,n\})$.
Let ${\cal K}_m(a_1,\dots,a_n)={\cal K}(a_{1,m},\dots,a_{n,m})\ (m\in\mathbb{N})$. Then there is $M\in\mathbb{N}$ such that $m>M$ implies that ${\cal K}_m(a_1,\dots,a_n)$ is mean-like.
\end{prp}
\begin{proof}Let $a_1=\min\{a_1,\dots,a_n\},\ a_n=\max\{a_1,\dots,a_n\},\ k={\cal K}(a_1,\dots,a_n)$.
${\cal K}$ is strict hence $a_1<k<a_n$. Choose $\epsilon>0$ such that $a_1<k-\epsilon<k+\epsilon<a_n$. Then by the continuity of ${\cal K}$ find $M\in\mathbb{N}$ such that $m>M$ implies that ${\cal K}(a_{1,m},\dots,a_{n,m})\in(k-\epsilon,k+\epsilon)$. 
\end{proof}

\begin{rem}The two variable $\min$ shows that if we omit the condition of strictness then the previous proposition fails to be true. 
\end{rem}

\subsection{A quasi-mean from physics}

The motivation of the next example comes from physics. Suppose we have an electrical network with two resistors connected into a parallel circuit. To find the total resistance of the components we have the following formula: 
\begin{equation}\label{etrans}
R_{total}=\frac{1}{\frac{1}{R_1}+\frac{1}{R_2}}
\end{equation}
where $R_1,R_2$ are the resistance of the two components.

\begin{ex}\label{ephelctr}For $a,b\in\mathbb{R}^+$ set
\[{\cal K}(a,b)=\frac{1}{\frac{1}{a}+\frac{1}{b}}.\]
This is a kind of perturbation of the harmonic mean for two variables.
\end{ex}

\begin{prp}Let ${\cal K}$ be the mean defined in \ref{ephelctr}. Then the following statements hold:

1. ${\cal K}(a,b)\leq\min\{a,b\}\leq\max\{a,b\}$ hence it is a strong right-mean

2. ${\cal K}$ is not reflexive (${\cal K}(a,a)=\frac{a}{2}$), it is symmetric, continuous, monotone

3. $\lim\limits_{x\to +\infty}{\cal K}(a,x)=a$

4. $\lim\limits_{x\to 0+0}{\cal K}(a,x)=0$. \qed
\end{prp}

Of course we can generalize this quasi-mean for several variables in two ways. First according to this physics example when applying the formula \ref{etrans} for several transistors (the second generalization comes in \ref{epowerm}): 
\[{\cal K}(a_1,\dots,a_n)=\frac{1}{\frac{1}{a_1}+\dots+\frac{1}{a_n}}\ \ (a_1,\dots,a_n\in\mathbb{R}^+).\]
Similarly to the two variable case it can be readily seen it is a strong right-mean, it is not reflexive and it is symmetric, continuous, monotone.

Here we mention one more property.

\begin{prp}Let a sequence $(a_n)$ be given such that $\varliminf a_n>0$. Then
\[\lim\limits_{n\to\infty}{\cal K}(a_1,\dots,a_n)=\lim\limits_{n\to\infty}\frac{1}{\sum\limits_{i=1}^n\frac{1}{a_i}}=0.\tag*{\qed}\]
\end{prp}

\subsection{Quasi-arithmetic quasi-means}\label{sqaqmeans}

\begin{df}
Let $n\in\mathbb{N}-\{1\},\ a_1,\dots,a_n\in\mathbb{R},\ f:\mathbb{R}\to\mathbb{R}$ be continuous and increasing. 
Let ${\cal L}$ be a right-mean (left-mean) defined on subsets of $\mathbb{R}^n$ for each $n\in\mathbb{N}-\{1\}$.
Set 
\[\text{Dom }{\cal K}_{f,{\cal L}}=\big\{(a_1,\dots,a_n)\in\mathbb{R}^n\ (n\in\mathbb{N}): (f(a_1),\dots,f(a_n))\in\text{Dom }{\cal L}\big\}.\]
\[{\cal K}_{f,{\cal L}}(a_1,\dots,a_n)=f^{-1}\big({\cal L}(f(a_1),\dots,f(a_n))\big)\ \ (m\in\mathbb{Z}).\]
\end{df}

\begin{prp}\label{pkflqm}If ${\cal L}$ is a right-mean (left-mean) then so does ${\cal K}_{f,{\cal L}}$.
\end{prp}
\begin{proof}We just prove the "right" part, the other is similar.

Let $a_j=\max\{a_1,\dots,a_n\}$. Then $f(a_j)=\max\{f(a_1),\dots,f(a_n)\}$ and ${\cal L}(f(a_1),\dots,f(a_n))\leq f(a_j)$ which gives that ${\cal K}_{f,{\cal L}}(a_1,\dots,a_n)\leq a_j$.
\end{proof}

\begin{prp}\label{pkflqm2}Let ${\cal L}$ be a quasi-mean. If ${\cal L}$ is strict, monotone, symmetric, continuous then so does ${\cal K}_{f,{\cal L}}$ respectively.

Moreover ${\cal L}$ is reflexive iff ${\cal K}_{f,{\cal L}}$ is reflexive.
\end{prp}
\begin{proof}All properties can be readily checked.
\end{proof}

We are going to mention two special cases.

First we can generalize the quasi-means \ref{ebessel} and \ref{ephelctr} easily.

\begin{ex}\label{eqaritm}Let $n\in\mathbb{N}-\{1\},\ a_1,\dots,a_n\in\mathbb{R}^+,\ f:\mathbb{R}^+\to\mathbb{R}^+$ be continuous and increasing. Set 
\[{\cal K}_{f,{\cal B}^+}(a_1,\dots,a_n)=f^{-1}\big({\cal B}^+(f(a_1),\dots,f(a_n))\big)=f^{-1}\left({\frac{f(a_1)+\dots +f(a_n)}{n-1}}\right).\]
\end{ex}

Especially starting from power means we get:

\begin{ex}\label{epowerm}Let $n\in\mathbb{N}-\{1\},\ a_1,\dots,a_n\in\mathbb{R}^+,\ x\in\mathbb{R}$. Set 
\[{\cal K}_x(a_1,\dots,a_n)=\sqrt[x]{\frac{a_1^x+\dots +a_n^x}{n-1}},\]
\[{\cal K}_0(a_1,\dots,a_n)=\sqrt[n-1]{a_1\cdots a_n}.\]
\end{ex}

\begin{prp}The following statements hold:

1.  $\min\{a_1,\dots,a_n\}\leq{\cal K}_{f,{\cal B}^+}(a_1,\dots,a_n)$ i.e. ${\cal K}_{f,{\cal B}^+}$ is a quasi-mean

2. ${\cal K}_{f,{\cal B}^+}$ is strict, monotone, symmetric, continuous but it is not reflexive.
\end{prp}
\begin{proof}${\cal B}^+$ is a left-mean and apply \ref{pkflqm} and \ref{pkflqm2}.
\end{proof}

\begin{rem}a) $\lim\limits_{x\to 0}{\cal K}_x={\cal K}_0$ does not hold moreover $\lim\limits_{x\to 0}{\cal K}_x(a_1,\dots,a_n)$ does not exist in general.

b) $x<y\not\Rightarrow{\cal K}_x\leq{\cal K}_y$ e.g. $\sqrt{2\cdot 3\cdot 4}\not\leq \frac{2+3+4}{2}\not\leq \sqrt{\frac{2^2+3^2+4^2}{2}}$.
\end{rem}

\begin{prp}If $n=2,0<v<u$ then ${\cal K}_u(a,b)<{\cal K}_v(a,b)\ (a,b\in\mathbb{R}^+)$.
\end{prp}
\begin{proof}Let $c=a^v,d=b^v,k=\frac{u}{v}$. Then we have to show that $(c+d)^k>c^k+d^k$. Change $c$ to $x$ and consider it as a variable, while let us fix $d$. Then when $x=0$ then equality holds. Derivating by $x$ we get that $k(x+d)^{k-1}>kx^{k-1}$ which obviously holds.
\end{proof}

For quasi-arithmetic means there is another way how one can define associated quasi-means in a similar manner than in the definition of \ref{dqmna}.

\begin{df}
Let $n\in\mathbb{N}-\{1\},\ a_1,\dots,a_n\in\mathbb{R},\ f:\mathbb{R}\to\mathbb{R}$ be continuous and increasing. Set 
\[\text{Dom }{\cal K}_{f,m}=\{(a_1,\dots,a_n)\in\mathbb{R}^n\ (n\in\mathbb{N}): \exists i\in\{1,\dots,n\}\text{ such that }\lfloor10^mf(a_i)\rfloor\ne 0\}.\]
\[{\cal K}_{f,m}(a_1,\dots,a_n)=f^{-1}\big(\A_m(f(a_1),\dots,f(a_n))\big)\ \ (m\in\mathbb{Z}).\]
\end{df}

\begin{prp}Let $n\in\mathbb{N}-\{1\},\ a_1,\dots,a_n\in\mathbb{R},\ f:\mathbb{R}\to\mathbb{R}$ be continuous and increasing. Then ${\cal K}_{f,m}$ is a right-mean that is strict, monotone, symmetric, continuous but not reflexive..
\end{prp}
\begin{proof}$\A_m$ is a right-mean and apply \ref{pkflqm} and \ref{pkflqm2}.
\end{proof}

\subsection{Various other quasi-means}

\begin{ex}Let an ordinary symmetric mean ${\cal M}$ is given that is defined on arbitrarily many positive numbers i.e. ${\cal M}:\mathbb{R}^n\to\mathbb{R}$ is defined for all $n\in\mathbb{N}$. Or we can also rephrase that as ${\cal M}$ is defined on all finite multi-subsets of $\mathbb{R}$.

Set ${\cal M}'(\{a_i:1\leq i\leq n\})={\cal M}(\{a_i>0:1\leq i\leq n\})$ and ${\cal M}'(\emptyset)=0$ i.e. we apply ${\cal M}$ on the positive elements only. Then e.g. ${\cal M}'(\{-1,-2\})=0$. Obviously $\min\{a_1,\dots,a_n\}\leq{\cal M}'(a_1,\dots,a_n)$ Hence it is a quasi-mean.

It can be readily seen that if ${\cal M}$ is continuous then ${\cal M}'$ is right-continuous. Left-continuity is still holds if $0$ is not among the elements. If it is then ${\cal M}'$ is not left-continuous.

Monotonicity does not hold either: ${\cal A}'(\{-1,2,3\})\not\leq{\cal A}'(\{1,2,3\})$.
\end{ex}

\begin{ex}Let ${\cal K}(a,b)=\frac{\sqrt{a^2+b^2}}{2}\ (a,b\in\mathbb{R}^+)$. Obviously ${\cal K}(a,b)\leq\max\{a,b\}$ hence it is a quasi-mean. It is strict, monotone, symmetric, continuous but it is not reflexive.

There is a smaller domain where ${\cal K}$ becomes a mean. One has to restrict the domain to the values $\big\{(a,b)\in(\mathbb{R}^+)^2:\sqrt{3}a\leq b\text{ or }\sqrt{3}b\leq a\big\}$. 
\end{ex}

\begin{ex}Let $a_1,\dots,a_n\in\mathbb{R}^+\cup\{0\}$. Then set
\[{\cal K}(a_1,\dots,a_n)=\frac{a_1+a_1a_2+a_1a_2a_3+\dots+a_1\cdots a_n}{n}\]
If $1\leq a_i (1\leq i\leq n)$ then $\min\{a_1,\dots,a_n\}\leq{\cal K}(a_1,\dots,a_n)$. Actually if $a_i=\min\{a_1,\dots,a_n\}$ then we get that 
\[a_i=\frac{na_i}{n}\leq\frac{a_i+a_i^2+\dots+a_i^n}{n}\leq{\cal K}(a_1,\dots,a_n).\]
Similarly we get that if $1\geq a_i\ (1\leq i\leq n)$ then ${\cal K}(a_1,\dots,a_n)\leq\max\{a_1,\dots,a_n\}$. Therefore ${\cal K}|_{[1,+\infty)}$ is a right-mean while ${\cal K}|_{[0,1]}$ is a left-mean. 

Clearly they are not means simply because they are not reflexive.
\end{ex}

\begin{ex}Let $a_1,\dots,a_n\in\mathbb{R}^+\cup\{0\}$. Then set
\[{\cal K}(a_1,\dots,a_n)=\sqrt[n]{\frac{a_1+a_1a_2+a_1a_2a_3+\dots+a_1\cdots a_n}{n}}\]
If $1\geq a_i (1\leq i\leq n)$ then $\min\{a_1,\dots,a_n\}\leq{\cal K}(a_1,\dots,a_n)$. Actually if $a_i=\min\{a_1,\dots,a_n\}$ then we get that 
\[a_i=\sqrt[n]{\frac{na_i^n}{n}}\leq\sqrt[n]{\frac{a_i+a_i^2+\dots+a_i^n}{n}}\leq{\cal K}(a_1,\dots,a_n).\]
Similarly we get that if $1\leq a_i (1\leq i\leq n)$ then ${\cal K}(a_1,\dots,a_n)\leq\max\{a_1,\dots,a_n\}$. Therefore ${\cal K}|_{[1,+\infty)}$ is a left-mean while ${\cal K}|_{[0,1]}$ is a right-mean. 

Clearly it is not a mean simply because they are not reflexive. 
\end{ex}

\begin{ex}a) Let ${\cal K}(a,b)=\frac{a+b}{2+\max\{a,b\}-\min\{a,b\}}\ (a,b\in\mathbb{R}^+)$. Obviously ${\cal K}(a,b)\leq\max\{a,b\}$ hence it is a quasi-mean. ${\cal K}(2,4)=1.5$ shows that it is not a mean. It is strict, symmetric, continuous and reflexive hence it is mean-continuous as well. It is not monotone as $\frac{2+3}{3}\not\leq\frac{2+4}{4}$ witnesses.

b) If $a,b\in\mathbb{R}^+$ let 
\[{\cal K}(a,b)=
\begin{cases}
\frac{a+b}{2+\frac{1}{\max\{a,b\}-\min\{a,b\}}}&\text{if }a\ne b\\
a&\text{otherwise.}
\end{cases}
\]
Obviously ${\cal K}(a,b)\leq\max\{a,b\}$ hence it is a quasi-mean. ${\cal K}(10,11)=7$ shows that it is not a mean. It is strict, symmetric, reflexive, not continuous and not mean-continuous.
\end{ex}

\begin{ex}It is easy to create a non-continuous, non-reflexive but mean-continuous quasi means. E.g. if $a,b\in[1,+\infty)$ let 
\[{\cal K}(a,b)=
\begin{cases}
\frac{a+b}{2}&\text{if }a\ne b\\
a+\frac{1}{a}&\text{otherwise.}
\end{cases}
\]
\end{ex}

\section{Duality and quasi-means}

There is a kind of algebraic duality between certain means that are defined by addition and multiplication only (or some other operation that is derived from those operations). If we have such a mean then if in the definition of the mean we replace addition with multiplication and vice versa we may end up with a new mean that is called the dual of the original mean. Unfortunately this operation is not well defined because the same mean can be represented in many different formulas. All of those can be investigated and see if the dual is a mean or not. In some cases this operation results quasi-means. 

\smallskip

First let us give some non proper explanation why we will replace $\frac{a}{n}$ with $\sqrt[n]{a}$. Let $x=\frac{a}{n}$. Then $nx=a$ that is $x+\dots + x=a$. Replacing $+$ with $\cdot$ we get $x^n$ hence $\frac{a}{n}\leadsto\sqrt[n]{a}$. Similar reasoning yields that $\sqrt[n]{a}\leadsto\frac{a}{n}$.

\begin{ex}If we consider the arithmetic mean and we replace addition with multiplication, division by $n$ with $n^{th}$ root then we end up with the geometric mean. And vice versa: from the geometric mean we get the arithmetic mean if we replace multiplication with addition, $n^{th}$ root with division by $n$:
\[\frac{a_1+\dots +a_n}{n}\leftrightsquigarrow\sqrt[n]{a_1\cdots a_n}.\]
\end{ex}

\begin{ex}If we start with the $x^{th}$ power mean ($x\in\mathbb{R}$) and we replace addition with multiplication, division by $n$ with $n^{th}$ root, power to $x$ with multiplication by $x$ then we end up with the geometric mean again:
\[\sqrt[x]{\frac{a_1^x+\dots +a_n^x}{n}}\leadsto\frac{\sqrt[n]{xa_1\cdots xa_n}}{x}=\sqrt[n]{a_1\cdots a_n}.\]
\end{ex}

We present an example which demonstrates that different but equivalent formulas can end up with different results.

\begin{ex}Let $a_1,a_2,a_3\in\mathbb{R}^+$.
\[\sqrt[6]{\frac{(a_1^2+a_2^2)(a_2^2+a_3^2)(a_3^2+a_1^2)}{8}}=\sqrt[3]{\sqrt{\frac{a_1^2+a_2^2}{2}}\sqrt{\frac{a_2^2+a_3^2}{2}}\sqrt{\frac{a_3^2+a_1^2}{2}}}\leadsto\]
\[\leadsto\frac{\sqrt{a_1a_2}+\sqrt{a_2a_3}+\sqrt{a_3a_1}}{3}\]
which is a mean. While
\[\sqrt[6]{\frac{(a_1^2+a_2^2)(a_2^2+a_3^2)(a_3^2+a_1^2)}{8}}\leadsto\frac{\sqrt[8]{(2a_12a_2)+(2a_22a_3)+(2a_32a_1)}}{6}=\]
\[=\frac{1}{2^{\frac{3}{4}}}\frac{a_1a_2+a_2a_3+a_3a_1}{3}\]
is only a quasi-mean.
\end{ex}

\section{Operations on quasi-means}

\begin{prp}Let ${\cal K}_0,\dots,{\cal K}_n$ be $n$-variable left-means (right-means). Set
\[{\cal K}(a_1,\dots,a_n)={\cal K}_0\big({\cal K}_1(a_1,\dots,a_n),\dots,{\cal K}_n(a_1,\dots,a_n)\big)\ \ (a_1,\dots,a_n\in\mathbb{R}).\]
Then ${\cal K}$ is a  left-mean (right-mean) as well.
\end{prp}
\begin{proof}We show the "right" part, the other is similar.

Let $a_n=\max\{a_1,\dots,a_n\}$ and let $j$ be chosen such that ${\cal K}_j(a_1,\dots,a_n)=\max\{{\cal K}_1(a_1,\dots,a_n),\dots,{\cal K}_n(a_1,\dots,a_n)\}$.
\[{\cal K}_0\big({\cal K}_1(a_1,\dots,a_n),\dots,{\cal K}_n(a_1,\dots,a_n)\big)\leq{\cal K}_j(a_1,\dots,a_n)\leq a_n.\qedhere\]

\end{proof}

We now present a generalization of compounding for quasi-means. We formulate it for right-means; a similar statement can be given for left-means.

\begin{thm}\label{tcompounding}Let ${\cal K},{\cal M}:(\mathbb{R}^+\cup\{0\})^2\to\mathbb{R}^+\cup\{0\}$ be $2$-variable right-means such that ${\cal K}\leq{\cal M}$. Moreover let ${\cal M}$ be strict and continuous. Let $a,b\in\mathbb{R}^+\cup\{0\},\ a\leq b$. Set $a_0=a,b_0=b,\ a_{n+1}={\cal K}(a_n,b_n),\ b_{n+1}={\cal M}(a_n,b_n)\ (n\in\mathbb{N})$. Then $(a_n),(b_n)$ converge to the same limit, say $d$ and with that ${\cal N}(a,b)=d$ becomes a strict right-mean.

Moreover if ${\cal K}$ is continuous as well then ${\cal N}$ is upper semi-continuous.
\end{thm}
\begin{proof}It can be readily seen that $\forall n\ a_n\leq b_n,\ b_{n+1}<b_n<b$. Hence $(b_n)$ is decreasing, bounded from below therefore converges, say to $d$. Let $c$ be an accumulation point of $(a_n)$. Let $a_{n_k}\to c$. Suppose that $c<d$. Let $e={\cal M}(c,d)<d$. Let $\epsilon>0$ be chosen such that $e+\epsilon<d$. Then by the continuity of ${\cal M}$ there is $K$ such that $n_k>K$ implies that ${\cal M}(a_{n_k},b_{n_k})\in(e-\epsilon,e+\epsilon)$ which would imply that $b_{n_k+1}<d$ that is a contradiction.

Clearly $d<b$ hence ${\cal N}(a,b)$ is a strict right-mean.

To see the upper semi-continuity of ${\cal N}$ first observe that for the above defined sequences for all $n\in\mathbb{N}$ both $a_n,b_n$ are continuous functions of $a,b$ i.e $a_n=a_n(a,b)$ can be considered as a function of $a,b$ that is a two-variable continuous function.

Now let $x,y\in\mathbb{R}^+\cup\{0\},\ x\leq y,\ x_n\to x,\ y_n\to y,\ d={\cal N}(x,y)$. Suppose indirectly that $d_n={\cal N}(x_n,y_n)\not\to d$. Then changing to a subsequence we can assume that ${\cal N}(x_n,y_n)\to c > d$. Let $\epsilon=\frac{c-d}{3}$. Then there is $N$ such that $b_N(x,y)<d+\epsilon$. By the continuity of $b_N$ at $(x,y)$ there is $M$ such that  $m>M$ implies that $b_N(x_m,y_m)<d+2\epsilon$. But $d_m<b_N(x_m,y_m)$ which is a contradiction.
\end{proof}

We provide a sufficient condition for ${\cal N}$ being continuous.

\begin{prp}Let $x,y\in\mathbb{R}^+\cup\{0\}$. With the notation and conditions of \ref{tcompounding} assume that $\forall\epsilon>0\ \exists\delta>0\ \exists M\in\mathbb{N}$ such that $m\geq M, |x-x'|<\delta, |y-y'|<\delta$ implies that $|a_m(x',y')-{\cal N}(x',y')|<\epsilon$. Then ${\cal N}$ is continuous at $(x,y)$.
\end{prp}
\begin{proof}We use the notations of the proof of \ref{tcompounding}. We have to prove the lower semi-continuity of ${\cal N}$.

Suppose indirectly that $d_n={\cal N}(x_n,y_n)\not\to d$. Then changing to a subsequence we can assume that ${\cal N}(x_n,y_n)\to c < d$. Let $\epsilon=\frac{c-d}{4}$. By condition choose $\delta>0, M\in\mathbb{N}$ such that $m\geq M, |x-x'|<\delta, |y-y'|<\delta$ implies that $|a_m(x',y')-{\cal N}(x',y')|<\epsilon$.
Then there is $N>M$ such that $a_N(x,y)>d-\epsilon$. By the continuity of $a_N$ at $(x,y)$ there is $N'$ such that  $k>N'$ implies that $|x-x_k|<\delta, |y-y_k|<\delta$ and $a_N(x_k,y_k)>d-2\epsilon$. But then $d_k>d-3\epsilon$ which is a contradiction.
\end{proof}

\begin{ex}We provide quasi-means that satisfy the conditions of \ref{tcompounding} but the compounding quasi-mean fails to be continuous.

Let ${\cal K}={\cal A}^-_0|_2,{\cal M}={\cal A}^-_1|_2$. Let $a=1,b=2$ and $c_n\to 1,0<c_n<1,d_n\to 2,1<d_n<2$. Let us denote the compounding mean by ${\cal N}$.

We state that $1\leq{\cal N}(1,2)$. Clearly $1\leq a_0,b_0$ (using the notation of \ref{tcompounding}). Suppose that $1\leq a_n,b_n$. Then $1\leq a_{n+1},b_{n+1}$. Hence we get that $1\leq\lim b_n$.

We show that $\lim\limits_{n\to\infty}{\cal N}(c_n,d_n)<1$ moreover $\forall n\ {\cal N}(c_n,d_n)<0.75$. If $n$ is fixed then denote the associated sequences by $\left(a_m^{(c_n,d_n)}\right),\left(b_m^{(c_n,d_n)}\right)$. One can readily get that $a_1^{(c_n,d_n)}=0.5,\ b_1^{(c_n,d_n)}\leq 1.5$; $a_2^{(c_n,d_n)}\leq 0.5,\ b_2^{(c_n,d_n)}\leq 1$; $a_3^{(c_n,d_n)}\leq 0.5,\ b_3^{(c_n,d_n)}\leq 0.75$. Therefore $\lim\limits_{m\to\infty} b_m^{(c_n,d_n)}<0.75$.
\end{ex}

\begin{df}Let ${\cal K}:H^2\to H$ be a symmetric right-mean. If $a,b\in H,a\leq b$ then let ${\cal K}^{(0)}(a,b)={\cal K}^{(0)}(b,a)={\cal K}(a,b),\ {\cal K}^{(n+1)}(a,b)={\cal K}^{(n+1)}(b,a)={\cal K}\big({\cal K}^{(n)}(a,b),b\big)\ (n\in\mathbb{N})$.
\end{df}

\begin{df}Let ${\cal K}:H^2\to H$ be a symmetric right-mean. ${\cal K}$ is called right-idempotent if ${\cal K}^{(1)}={\cal K}$.
\end{df}

Of course similar notions can be defined for the "left" version.

\begin{prp}Let ${\cal K}:H^2\to H$ be a symmetric right-mean. Then the following statements hold:
\begin{enumerate}\setlength\itemsep{0em}
\item For all $n\in\mathbb{N}\ \ {\cal K}^{(n)}$ is a right-mean too.
\item If ${\cal K}$ is right-idempotent then $\forall n\ {\cal K}^{(n)}={\cal K}$.
\item If $\forall a\in H\ \forall b\in H\ \lim\limits_{n\to\infty}{\cal K}^{(n)}(a,b)$ exists and finite then let us denote it by $\tilde{\cal K}(a,b)$. If ${\cal K}$ is continuous then $\tilde{\cal K}$ is a right-idempotent right-mean.
\end{enumerate}
\end{prp}
\begin{proof}(1) and (2) are trivial.

(3): $\tilde{\cal K}$ is obviously a right-mean since more generally the limit of right-means is a right-mean. 

If $a,b\in H$ then set $y=\tilde{\cal K}(a,b)$. By definition ${\cal K}^{(n)}(a,b)\to y$ which gives that ${\cal K}\big({\cal K}^{(n)}(a,b),b\big)\to {\cal K}(y,b)$ since ${\cal K}$ is continuous. But ${\cal K}\big({\cal K}^{(n)}(a,b),b\big)={\cal K}^{(n+1)}(a,b)\to y$ too thus ${\cal K}(y,b)=y$. Which yields that $\tilde{\cal K}(y,b)=y$ too. Then we get that $\tilde{\cal K}\big(\tilde{\cal K}(a,b),b\big)=\tilde{\cal K}\big(y,b\big)=y=\tilde{\cal K}(a,b)$.
\end{proof}

\begin{ex}${\cal K}$ being continuous does not imply that $\tilde{\cal K}$ is continuous.

E.g. let ${\cal K}(a,b)=\min\{a^2,b\}\ \ (a,b\in[0,2])$. Then $\tilde{\cal K}(1,2)=1$ and $\tilde{\cal K}(1-\epsilon,2)=0$ for all $\epsilon>0$.
\end{ex}

In \cite{laesv} we investigated how one can extend a $k$-variable mean to an $n$-variable mean for $2\leq k<n$. Here we are going to generalize partially that process for quasi-means. We formulate the following theorems for right-means; one can easily get them for left-means as well.

\begin{prp}\label{pext1}Let $H\subset\mathbb{R}$ that is bounded from below. Let ${\cal K}:H^2\to H$ be a strict, continuous, monotone right-mean. Let $a,b,c\in H,\ a\leq b\leq c$. Let $a_0=a,b_0=b,c_0=c,\ a_{n+1}={\cal K}(a_n,b_n),\ b_{n+1}={\cal K}(a_n,c_n),\ c_{n+1}={\cal K}(b_n,c_n)\ (n\in\mathbb{N})$. Now all three sequences converge to the same limit.
\end{prp}
\begin{proof}First note that monotonicity gives that $\forall n\ a_n\leq b_n\leq c_n$. ${\cal K}$ is a right-mean which yields that $(c_n)$ is decreasing hence it is convergent by boundedness; let $c_n\to d$ .

Assume that $b_{n_k}\to e<d$. Then $c_{n_k+1}={\cal K}(b_{n_k},c_{n_k})$ but continuity yields that ${\cal K}(b_{n_k},c_{n_k})\to {\cal K}(e,d)<d$ by strictness which is a contradiction.

Similarly assume that $a_{n_k}\to e<d$. Then $b_{n_k+1}={\cal K}(a_{n_k},c_{n_k})$ but continuity yields that ${\cal K}(a_{n_k},c_{n_k})\to {\cal K}(e,d)<d$ by strictness which is a contradiction.
\end{proof}

\begin{ex}If one omits the boundedness property then \ref{pext1} fails to be true.

E.g. let ${\cal K}(a,b)=\frac{a+b}{2}-1$. It can be readily seen that $c_n\to -\infty$.
\end{ex}

\begin{prp}Let ${\cal K}:H^2\to H$ be a strict, continuous, monotone right-mean. With the notations of \ref{pext1} if $c_n\to -\infty$ then $a_n\to -\infty$ and $b_n\to -\infty$ as well. \qed
\end{prp}

\begin{prp}Let ${\cal K}:H^2\to H$ be a strict, continuous, monotone, semi-reflexive right-mean. With the notations of \ref{pext1} all three sequences converge to the same limit.
\end{prp}
\begin{proof}By induction we show that ${\cal K}(a,a)\leq a_n,b_n,c_n$ that is enough referring to the proof of \ref{pext1}.

Monotonicity gives that ${\cal K}(a,a)\leq a_1,b_1,c_1$. Then suppose that ${\cal K}(a,a)\leq a_n,b_n,c_n$. Then ${\cal K}(a,a)={\cal K}\big({\cal K}(a,a),{\cal K}(a,a)\big)\leq{\cal K}(a_n,b_n),{\cal K}(a_n,c_n),{\cal K}(b_n,c_n)$.
\end{proof}

\begin{thm}\label{text1}Let ${\cal K}:H^2\to H$ be a strict, continuous, monotone right-mean that is either semi-reflexive or $H$ is bounded from below. With the notations of \ref{pext1} let $\tilde{\cal K}(a,b,c)$ denote the common limit point of the sequences. Then $\tilde{\cal K}$ is a strict, upper semi-continuous, monotone right-mean.
\end{thm}
\begin{proof}$\tilde{\cal K}$ is a right-mean because $\tilde{\cal K}(a,b,c)\leq c_1\leq c$.

$\tilde{\cal K}$ is a strict since $\lim c_n\leq c_1<c$.

$\tilde{\cal K}$ is a monotone: Let $a'\leq a,b'\leq b,c'\leq c$ and the associated sequences are $(a'_n),(b'_n),(c'_n)$. Straightforward induction yields that $a'_n\leq a_n,b'_n\leq b_n,c'_n\leq c_n$ which gives that $\lim c'_n\leq \lim c_n$.

$\tilde{\cal K}$ is upper semi-continuous: First note that for a fixed $n\in\mathbb{N}$ $c_n$ can be considered as a function of $a,b,c$. Clearly $c_n(a,b,c)$ is a continuous 3-variable function. Now if $\varepsilon>0$ is given then $\exists n\in\mathbb{N}$ such that $c_n<\tilde{\cal K}(a,b,c)+\frac{\varepsilon}{2}$ ($c_n=c_n(a,b,c)$). And $\exists\delta>0$ such that $|a'-a|<\delta,|b'-b|<\delta,|c'-c|<\delta$ implies that $|c_n(a',b',c')-c_n|<\frac{\varepsilon}{2}$. Therefore $c_n(a',b',c')<\tilde{\cal K}(a,b,c)+\varepsilon$. We know that $\tilde{\cal K}(a',b',c')\leq c_n(a',b',c')$ holds and $\varepsilon$ being arbitrary gives that $\tilde{\cal K}(a',b',c')\leq\tilde{\cal K}(a,b,c)$.
\end{proof}

\begin{ex}\label{eA0c}Let us count what we get for $\A^-_0$ and $1.1,2.1,3.1$ if we apply the method described in \ref{pext1}.

1.step: 1.5,2,2.5;
2.step: 1.5,1.5,2;
3.step: 1,1.5,1.5;
4.step: 1,1,1.

Hence this extension does not give back $\A^-_0(1.1,2.1,3.1)=2$ in contrast with \A\ where it does give back the 3-variable version (see \cite{laesv}). 
\end{ex}

\begin{ex}$\A^-_0$ is a quasi-mean that satisfies the conditions of \ref{text1} but the extension quasi-mean $\tilde{\cal A}^-_0$ fails to be continuous. 

Let $a=1,b=2,c=3$. Then easy calculation shows that $\tilde{\cal A}^-_0(1,2,3)=1$.

Let $(a_n),(b_n),(c_n)$ three sequences such that $a_n\to 1,0<a_n<1,b_n\to 2,1<b_n<2,c_n\to 3,2<c_n<3$. Let us calculate $\tilde{\cal A}^-_0(a_n,b_n,c_n)$ for a fixed $n$.

1.step: 0.5,1,1.5;
2.step: 0.5,0.5,1;
3.step: 0,0.5,0.5;
4.step: 0,0,0.

Hence $\forall n\ \tilde{\cal A}^-_0(a_n,b_n,c_n)=0$ and $\lim\limits_{n\to\infty}\tilde{\cal A}^-_0(a_n,b_n,c_n)=0\ne \tilde{\cal A}^-_0(1,2,3)$.
\end{ex}

The phenomenon in the example \ref{eA0c} is not an accident. For $\A^-_m,\A^+_m$ these sequences become constant after finitely many steps as the next statement points out.

\begin{prp}Let ${\cal K}:H^2\to H$ be a strict, continuous, monotone, semi-reflexive right-mean such that $\mathrm{Ran\ }{\cal K}$ is closed, discrete. With the notations of \ref{pext1} $\exists N\in\mathbb{N}\ \exists d\in\mathbb{R}$ such that $n\geq N$ implies that $a_n=b_n=c_n=d$.
\end{prp}
\begin{proof}If for an $n\in\mathbb{N}$ not all values $a_n,b_n,c_n$ are equal then in the next step one of the elements $a_{n+1},b_{n+1},c_{n+1}$ will strictly decrease by the strictness of ${\cal K}$. Simple  induction gives that ${\cal K}(a,a)\leq a_n,b_n,c_n$ by monotonicity and semi-reflexivity. $\mathrm{Ran\ }{\cal K}$ being closed and discrete yields that there are finitely many values of $\mathrm{Ran\ }{\cal K}$ between ${\cal K}(a,a)$ and $c$ from which we get the claim.
\end{proof}

\section{On quasi-monotonicity}

We are going to generalize the usual notion of monotonicity.

\begin{df}A two variable right-mean ${\cal K}:H^2\to H$ is called \textbf{quasi-monotone} if $\forall a\in H\ \forall b\in H, a\leq b$ then the following sequence $(a_n)$ is monotone: Let $a_0=a,\ a_{n+1}={\cal K}(a_n,b)$.

A two variable left-mean ${\cal K}:H^2\to H$ is quasi-monotone if $\forall a\in H\ \forall b\in H, a\leq b$ then the following sequence $(b_n)$ is monotone: Let $b_0=b,\ b_{n+1}={\cal K}(a,b_n)$.
\end{df}

\begin{prp}If ${\cal K}:H^2\to H$ is a mean then it is quasi-monotone.
\end{prp}
\begin{proof}A mean is a right-mean and left-mean at the same time thus we have to check both regarding quasi-monotonicity. We will show the "right" part, the other is similar.

If $a,b\in H,\ a\leq b$ then $\forall n\in\mathbb{N}\ a_n\leq{\cal K}(a_n,b)=a_{n+1}\leq b$ hence $(a_n)$ is monotone.
\end{proof}

\begin{prp}If a quasi-mean ${\cal K}:H^2\to H$ is monotone then it is quasi-monotone as well.
\end{prp}
\begin{proof}We show it for right-means, the other is similar.

Suppose that $a_1\leq a_0$. Then $a_2={\cal K}(a_1,b)\leq{\cal K}(a_0,b)=a_1$ by monotonicity. And by induction one gets that $(a_n)$ is decreasing.

When $a_0\leq a_1$ then a similar argument shows that $(a_n)$ is increasing.
\end{proof}

\begin{ex}The converse is not true. 

Let ${\cal K}:[0,1]^2\to [0,1]$ be symmetric and defined as follows. If $a,b\in[0,1],\ a\leq b$ then set
\[{\cal K}(a,b)=
\begin{cases}
\frac{a}{2}&\text{if }b=1\text{ and }a=\frac{1}{2^n}\text{ for some }n\in\mathbb{N}\\
b&\text{otherwise.}
\end{cases}
\]
It can be readily seen that ${\cal K}$ is a quasi-monotone right-mean. It is not monotone as $1={\cal K}(0.4,1)\not\leq{\cal K}(0.5,1)=0.25$ demonstrates.
\end{ex}

\begin{ex}\label{ecsqm}The previous example was not continuous. Now we give a continuous counter example too.

Let ${\cal K}:[0,1]^2\to [0,1]$ be symmetric and defined as follows. If $a,b\in[0,1],\ a\leq b$ then set
\[{\cal K}(a,b)=
\begin{cases}
b\cdot 2\left(\frac{a}{b}\right)^2&\text{if }a\leq\frac{b}{2}\\
b\left(1-2\left(1-\frac{a}{b}\right)^2\right)&\text{otherwise.}
\end{cases}
\]
It can be readily checked that ${\cal K}$ is a continuous, quasi-monotone right-mean. It is not monotone as $0,2={\cal K}(0.2,0,4)\not\leq{\cal K}(0.2,1)=0.08$ demonstrates.
\end{ex}

\begin{df}A quasi-mean ${\cal K}:H^2\to H$ is called \textbf{right-injective} if $\forall b\in H\ {\cal K}(x,b)={\cal K}(y,b)$ implies that $x=y$. 

It is left-injective if $\forall a\in H\ {\cal K}(a,x)={\cal K}(a,y)$ implies that $x=y$.
\end{df}

\begin{ex}The example in \ref{ecsqm} is clearly right-injective but not left-injective as e.g. ${\cal K}(0,0.4)={\cal K}(0,1)=0$ shows.
\end{ex}

\begin{prp}A quasi-mean ${\cal K}:H^2\to H$ is right-injective iff $\forall b\in H\ f_b^{\cal K}(x)={\cal K}(x,b)\ \ (0\leq x\leq b)$ is an injective function.\qed
\end{prp}

\begin{prp}\label{pscqmrm}Let ${\cal K}:[0,1]^2\to [0,1]$ be a continuous, quasi-monotone right-mean. Set $f_b(x)={\cal K}(x,b)\ \ (0\leq x\leq b,\ b\in[0,1])$. Then the following statements hold:
\begin{enumerate}\setlength\itemsep{0em}
\item $Z_b=\{x\in[0,b]:f_b(x)=x\}$ is a closed set.
\item Let us use the notation $[0,b]-Z_b=\sideset{}{^*}\bigcup\limits_{n=1}^{\infty}N_i$ where $N_i=(a_i,b_i)$. Then either $\forall x\in N_i\ f_b(x)>x$ or $\forall x\in N_i\ f_b(x)<x$. Set $s(i)=1$ in the first case, $s(i)=-1$ in the second.
\item If $f_b(N_i)\cap N_j\ne\emptyset$ then $s(i)=s(j)$.
\item If moreover ${\cal K}$ is right-injective then $\forall i\ f_b|_{N_i}$ is a monotone function and $f_b(N_i)=N_i$. 
\end{enumerate}
\end{prp}
\begin{proof}
\begin{enumerate}\setlength\itemsep{0em}
\item Obvious because of the continuity of $f_b$.
\item The function $f_b(x)-x$ is continuous and never 0 on $N_i$.
\item If $\sup N_i\leq\inf N_j$ and $f_b(N_i)\cap N_j\ne\emptyset$ then $s(i)=1$ and there is $x\in N_i$ such that $f_b(x)\in N_j$. Then we get that $x<f_b(x)$ hence by the quasi-monotonicity $f_b(x)<f_b(f_b(x))$ which gives that $s(j)=1$ too.

The other case is similar.
\item Obvious. \qedhere
\end{enumerate}
\end{proof}

\begin{ex}Using the notation of \ref{pscqmrm} $b\mapsto Z_b$ is not continuous in the Hausdorff-metric in general.

E.g. let ${\cal K}(a,b)=a-ab+a^2b\ (a,b\in[0,1])$. Then clearly $Z_0=[0,1]$ and $Z_b=\{0,1\}$ if $b>0$.
\end{ex}

\section{Further generalizations}

We can go even further on generalization. Let us consider for example the following function.

\begin{df}Set $\A^*_m=\frac{\A^+_m+\A^-_m}{2}\ \ (m\in\mathbb{Z})$.
\end{df}

Obvious calculation shows that \[\A^*_m(2\cdot 10^{-m},2.1\cdot 10^{-m})=\frac{2.5\cdot 10^{-m}+2\cdot 10^{-m}}{2}=2.25\cdot 10^{-m}>\max\{2\cdot 10^{-m},2.1\cdot 10^{-m}\},\]
\[A^*_m(1.9\cdot 10^{-m},2\cdot 10^{-m})=\frac{2\cdot 10^{-m}+1.5\cdot 10^{-m}}{2}=1.75\cdot 10^{-m}<\min\{1.9\cdot 10^{-m},2\cdot 10^{-m}\}.\] 
Therefore $\A^*_m$ is not a quasi-mean $(m\in\mathbb{Z})$. However as $m$ increases it gets closer and closer to being a mean in some sense. In that direction we introduce the following terminology.

\begin{df}${\cal K}:H^n\to\mathbb{R}$ is called an \textbf{a-quasi-mean} if $\exists K\geq 0$ such that ${\cal K}(a_1,\dots,a_n)\in\big[\min\{a_1,\dots,a_n\}-K,\ \max\{a_1,\dots,a_n\}+K\big]$.

If ${\cal K}$ is a multi-variable quasi-mean then it is an a-quasi-mean if $\exists K\geq 0$ such that for all $n$ the above condition holds.\qed
\end{df}

\begin{prp}$\A^*_m$ is an a-quasi-mean with $K=\frac{1}{2\cdot 10^m}$.
\end{prp}
\begin{proof}By \ref{pqmna1} (5) and \ref{pqmna12} (5)  $\frac{\A^+_m+\A^-_m}{2}\leq\frac{\A+\frac{1}{10^m}+\A}{2}=\A+\frac{1}{2\cdot 10^m}\leq\max+\frac{1}{2\cdot 10^m}$. The other inequality is similar.
\end{proof}

\begin{prp}$\lim\limits_{m\to\infty}{\cal A}^*_m(a_1,\dots,a_n)={\cal A}(a_1,\dots,a_n)$ and the convergence is uniform ($(a_1,\dots,a_n)\in\text{Dom }\A^+_0\cap\text{Dom }\A^-_0$)

\end{prp}
\begin{proof}\ref{pqmna1} (6) and \ref{pqmna12} (6).
\end{proof}

\begin{prp}${\cal K}$ is an a-quasi-mean iff $\mathrm{mdist}({\cal K})<\infty$.\qed
\end{prp}

Now let us examine another quasi-mean like function. It is related to ${\cal B}^+$ and ${\cal B}^-$ with the difference that we do not restrict its domain.

\begin{df}Set ${\cal B}(a_1,\dots,a_n)=\frac{\sum\limits_{i=1}^na_i}{n-1}\ \ (n\in\mathbb{N},\ n\geq 2,\ a_i\in\mathbb{R})$.
\end{df}

Clearly it is not a quasi-mean as $2<{\cal B}(1,2)$ and ${\cal B}(-1,-2)<-1$ show.

\begin{df}${\cal K}:H^n\to\mathbb{R}$ is called an \textbf{m-quasi-mean} if $\exists K\geq 0$ such that 
\[{\cal K}(a_1,\dots,a_n)\in\big[\min\{a_1,\dots,a_n\}-K\max\{|a_1|,\dots,|a_n|\},\]
\[\hspace{3.3cm}\max\{a_1,\dots,a_n\}+K\max\{|a_1|,\dots,|a_n|\}\big].\]

If ${\cal K}$ is a multi-variable quasi-mean then it is an m-quasi-mean if $\exists K\geq 0$ such that for all $n$ the above condition holds.\qed
\end{df}

\begin{prp}${\cal B}|_n$ is an m-quasi-mean with $K=\frac{1}{n-1}$.
\end{prp}
\begin{proof}Let $a_n=\max\{a_1,\dots,a_n\}$ and $0\leq a_n$. Then 
\[{\cal B}(a_1,\dots,a_n)\leq\frac{n}{n-1}a_n=a_n+\frac{1}{n-1}a_n\leq a_n+\frac{1}{n-1}\max\{|a_1|,\dots,|a_n|\}.\]

If $a_n<0$ then $a_n+\frac{1}{n-1}a_n\leq a_n$.

Similar arguments work for the other inequality.
\end{proof}

\begin{cor}${\cal B}$ is an m-quasi-mean with $K=1$.\qed
\end{cor}

\begin{prp}If $\mathrm{mdistp}({\cal K})<\infty$ then ${\cal K}$ is an m-quasi-mean.
\end{prp}
\begin{proof}Let $(a_1,\dots,a_n)\in\mathrm{Dom\ }{\cal K}$. Set $m=\min\{a_1,\dots,a_n\},\ M=\max\{a_1,\dots,a_n\},\ N=\max\{|a_1|,\dots,|a_n|\}$.
Suppose that $\exists K$ such that ${\cal K}(a_1,\dots,a_n)-M\leq K(M-m)$. Then 
\[{\cal K}(a_1,\dots,a_n)\leq M+KM-Km\leq M+KN+KN=M+(K+1)N.\qedhere\]
\end{proof}

Obvious examples show that a quasi-mean is not necessarily an a-quasi-mean or an m-quasi-mean. And vice versa an a-quasi-mean or an m-quasi-mean is not necessarily a quasi-mean.

\begin{ex}An a-quasi-mean is not necessarily an m-quasi-mean, e.g. let ${\cal K}(a,b)=\max\{a,b\}+1\ \ (a,b\in\mathbb{R})$. If $a=0,b>0$ then $K=\frac{1}{b}$ shows the it is not an m-quasi-mean.

\smallskip

An m-quasi-mean is not necessarily an a-quasi-mean, e.g. let ${\cal K}(a,b)=2\cdot\max\{a,b\}\ \ (a,b\geq 1)$.
\end{ex}


{\footnotesize
\noindent
Dennis G\'abor College, Hungary 1119 Budapest Fej\'er Lip\'ot u. 70.

\noindent E-mail: losonczi@gdf.hu, alosonczi1@gmail.com\\
}
\end{document}